\documentclass[12pt]{amsart}
\def\to{\rightarrow}
\newtheorem{theorem} {Theorem}
\newtheorem{corollary} {Corollary}
\theoremstyle{plain}
\begin{document}
\title [Partial sums]
{Partial sums of hypergeometric series of unit argument}
\author 
{Wolfgang B\"uhring}
\address 
{Physikalisches Institut\\ Universit\"at Heidelberg\\ Philosophenweg~12\\ 69120~Heidelberg\\
Germany}
\email 
{buehring@physi.uni-heidelberg.de}
\abstract 
The asymptotic behaviour of partial sums of generalized hypergeometric series of unit argument is investigated.  
\endabstract
\maketitle
\section{Introduction}
This paper deals with finite sums
\begin{equation}\sum\limits_{l=0}^{m-1} {{{\Gamma (a_1+l)\Gamma (a_2+l)\cdots \Gamma (a_{p+1}+l)} \over {\Gamma (b_1+l)\cdots \Gamma (b_p+l)\Gamma (1+l)}}}
\label{e1}
\end{equation}
which are partial sums of the first $m$ terms of hypergeometric series ${}_{p+1}F_p$ of unit argument $z=1$ multiplied by appropriate gamma function factors,
\begin{equation}
{{\Gamma (a_1)\Gamma (a_2)\cdots \Gamma (a_{p+1})} \over {\Gamma (b_1)\cdots \Gamma (b_p)}}{}_{p+1}F_p\left( \left.\begin{array}{c}a_1,a_2,\ldots ,a_{p+1}\\
  b_1,\ldots,b_p\end{array} \right| z \right)
\label{e2}
\end{equation}
\[=\sum\limits_{l=0}^\infty  {{{\Gamma (a_1+l)\Gamma (a_2+l)\cdots \Gamma (a_{p+1}+l)} \over {\Gamma (b_1+l)\cdots \Gamma (b_p+l)\Gamma (1+l)}}z^l}.\]
An important quantity in this context is 
\begin{equation}
s_p=b_1+\cdots +b_p-a_1-a_2-\cdots -a_{p+1},
\label{e3}
\end{equation}
the non-trivial characteristic exponent of the underlying differential equation at $z=1$.

If $\Re(s_p)>0$, then the hypergeometric series converges at $z=1$, and its value is given by the Gaussian summation formula for $p=1$ or by the corresponding generalized formulas \cite{b92} for $p>1$. The terms of the series behave asymptotically as
 \[{l^{a_1+a_2+\cdots +a_{p+1}-b_1-\cdots -b_p-1}}=l^{-s_p-1}\]
 for $l\to\infty$, and if we consider the partial sum of the first $m$ terms, then the contribution of the missing tail is $O(m^{-s_p})$ as $m\to\infty$. Therefore we have
\begin{theorem} 
If  $\Re(s_1)>0$, then
\begin{equation}\sum\limits_{l=0}^{m-1} {{{\Gamma (a_1+l)\Gamma (a_2+l
)}} \over {\Gamma (b_1+l)\Gamma (1+l)}}={{\Gamma (a_1)\Gamma (a_2)\Gamma (s_1)} \over {\Gamma (a_1+s_1)\Gamma (a_2+s_1)}}+O(m^{-s_1})
\label{e4}
\end{equation}
as $m\to\infty$,
\label{th1}
\end{theorem}
or more generally,
\begin{theorem}  
If $\Re(s_p)>0$, then
\begin{equation}\sum\limits_{l=0}^{m-1} {{{\Gamma (a_1+l)\Gamma (a_2+l)\cdots \Gamma (a_{p+1}+l)} \over {\Gamma (b_1+l)\cdots \Gamma (b_p+l)\Gamma (1+l)}}}
=g_0(0)+O(m^{-s_p})
\label{e5}
\end{equation}
as $m\to\infty$, where
\begin{equation}
g_0(0)
={{\Gamma (a_1)\Gamma (a_2)\Gamma (s_p)} \over {\Gamma (a_1+s_p)\Gamma (a_2+s_p)}}\sum\limits_{k=0}^\infty  {{{(s_p)_k} \over {(a_1+s_p)_k(a_2+s_p)_k}}A_k^{(p)}}
\label{e6}
\end{equation}
and the sum converges provided that 
 $\Re(a_j)>0$ for $j=3, \ldots , p+1$.
\label{th2}
\end{theorem}
Here use is made of the Pochhammer symbol
\[ (x)_n=x(x+1)\cdots (x+n-1)=\Gamma (x+n)/\Gamma (x), \]
and the coefficients $A_k^{(p)}$ for $p=2, 3, \ldots$ are given in \cite{b92}, but a few of them are displayed here again for convenience:
\begin{equation}
A_k^{(2)}={{(b_2-a_3)_k(b_1-a_3)_k} \over {k! }},
\label{e7}
\end{equation}
\begin{equation}
A_k^{(3)}=\sum\limits_{k_2=0}^k {{{(b_3+b_2-a_4-a_3+k_2)_{k-k_2}(b_1-a_3)_{k-k_2}(b_3-a_4)_{k_2}(b_2-a_4)_{k_2}} \over {(k-k_2)! k_2!}}},
\label{e8}
\end{equation}
\begin{equation}
A_k^{(4)}=\sum\limits_{k_2=0}^k {{{(b_4+b_3+b_2-a_5-a_4-a_3+k_2)_{k-k_2}(b_1-a_3)_{k-k_2}} \over {(k-k_2)! }}}
\label{e9}
\end{equation}
  \[\times \sum\limits_{k_3=0}^{k_2} {{{(b_4+b_3-a_5-a_4+k_3)_{k_2-k_3}(b_2-a_4)_{k_2-k_3}} \over {(k_2-k_3)! }}}\]
  \[\times {{(b_4-a_5)_{k_3}(b_3-a_5)_{k_3}} \over {k_3! }}.\]
 For $p=3, 4, ...$ several other representations are possible \cite{b92}, such as
\begin{equation}
A_k^{(3)}={{(b_3+b_2-a_4-a_3)_k(b_1-a_3)_k} \over {k! }}
\label{e10}
\end{equation}
\[\times \,_3F_2\left( \left. \begin{array}{c}b_3-a_4,b_2-a_4,-k\\
  b_3+b_2-a_4-a_3,1+a_3-b_1-k\end{array} \right|1 \right)\]
or
\begin{equation}
A_k^{(3)}={{(b_1+b_3-a_3-a_4)_k(b_2+b_3-a_3-a_4)_k} \over {k! }}
\label{e11}
\end{equation}
\[\times \,_3F_2\left( \left. \begin{array}{c}b_3-a_3,b_3-a_4,-k\\
  b_1+b_3-a_3-a_4,b_2+b_3-a_3-a_4\end{array} \right|1 \right) \]
and
\begin{equation}
A_k^{(4)}={{(b_4+b_3+b_2-a_5-a_4-a_3)_k(b_1-a_3)_k} \over {k! }}
\label{e49}
\end{equation}
\[\times \sum\limits_{l=0}^k {{{(b_4+b_3-a_5-a_4)_l(b_2-a_4)_l(-k)_l} \over {(b_4+b_3+b_2-a_5-a_4-a_3)_l(1+a_3-b_1-k)_ll! }}}\]
\[\times \,_3F_2\left( \left. \begin{array}{c}b_4-a_5,b_3-a_5,-l\\
  b_4+b_3-a_5-a_4,1+a_4-b_2-l\end{array} \right|1 \right) \]
or
\begin{equation}
A_k^{(4)}={{(b_1+b_3+b_4-a_3-a_4-a_5)_k(b_2+b_3+b_4-a_3-a_4-a_5)_k} \over {k! }}
\label{e50}
\end{equation}
\[\times \sum\limits_{l=0}^k {{{(b_3+b_4-a_3-a_5)_l(b_3+b_4-a_4-a_5)_l(-k)_l} \over {(b_1+b_3+b_4-a_3-a_4-a_5)_l(b_2+b_3+b_4-a_3-a_4-a_5)_ll! }}}\]
\[\times \,_3F_2\left( \left. \begin{array}{c}b_3-a_5,b_4-a_5,-l\\
  b_3+b_4-a_3-a_5,b_3+b_4-a_4-a_5\end{array} \right|1 \right). \]

In the case $p=2$, the sum in (\ref{e6}) is an ${}_3F_2$ and so (\ref{e5})--(\ref{e6}) reduce, when $m \to\infty$, to a well-known transformation formula for ${}_3F_2$ of unit argument \cite{bai} \cite{luk} \cite{sla}.

The simple formula (\ref{e4}) above corresponding to $p=1$ can be recovered from the general formula (\ref{e5})--(\ref{e6}) if we define 
\begin{equation}
A_0^{(1)}=1, \quad A_k^{(1)}=0 \quad \mbox{for} \quad  k>0,
\end{equation}
so that then the sum over $k$ is equal to $1$ and disappears. In a similar way the formulas in the other theorems below simplify for $p=1$.

It is the purpose of this work to investigate the asymptotic behaviour of the partial sums (\ref{e1}) for situations when $\Re(s_p)\le 0$ and to find a more detailed formula in cases, like  $s_p=1$, when a few additional terms of higher order might be desirable.

\section{The case when $s_p$ is not equal to an integer}

If $\Re(s_p)\le 0$, then the hypergeometric series diverges at $z=1$, and the contribution of the missing tail of the series does not asymptotically vanish. Responsible for this behaviour are the singular terms of the hypergeometric function at $z=1$. We have to explicitly subtract  the contribution from the most important singular terms in order to get the tail to vanish as $m\to\infty$. The behaviour of the hypergeometric function at $z=1$ is given by the continuation formula
\begin{equation}
{{\Gamma (a_1)\Gamma (a_2)\cdots \Gamma (a_{p+1})} \over {\Gamma (b_1)\cdots \Gamma (b_p)}}{}_{p+1}F_p\left( \left.\begin{array}{c}a_1,a_2,\ldots ,a_{p+1}\\
  b_1,\ldots,b_p\end{array} \right| z \right)
\label{e12}
\end{equation} 
\[=\sum\limits_{n=0}^\infty  {g_n(0)\,(1-z)^n}+(1-z)^{s_p}\sum\limits_{n=0}^\infty  {g_n(s_p)(1-z)^n},
\]
which is valid provided that $s_p$ is not equal to an integer. Here the coefficient $g_0(0)$ defined in (\ref{e6}) above enters, and, while the other 
$g_n(0)$ are not needed here, the coefficients $g_n(s_p)$ of the singular term are \cite{b92}
\begin{equation}
g_n(s_p)=(-1)^n{{(a_1+s_p)_n(a_2+s_p)_n\Gamma (-s_p-n)} \over {(1)_n}}
\label{e13}
\end{equation}
\[\times \sum\limits_{k=0}^n  {{{(-n)_k} \over {(a_1+s_p)_k(a_2+s_p)_k}}A_k^{(p)}}.\]
The singular terms on the right-hand side of (\ref{e12}) have, with $x=s_p+n$, the $z$-dependence
\begin{equation} (1-z)^x=\sum\limits_{l=0}^\infty  {{{ (-x)_l} \over {\Gamma (1+l)}}z^l}.
\label{e14}
\end{equation}
With $z=1$, we shall need the partial sum of this series,
\begin{equation}
\sum\limits_{l=0}^{m-1} {{{(-x)_l} \over {\Gamma (1+l)}}}=-{1 \over x}{{(-x)_m} \over {\Gamma (m)}}\quad (m=1,2,\ldots \quad x\ne 0).
\label{e15}
\end{equation}
Treating the first $N+1$ singular terms on the right-hand side of (\ref{e12}) in this way, we may obtain
\begin{equation}\sum\limits_{l=0}^{m-1} {{{\Gamma (a_1+l)\Gamma (a_2+l)\cdots \Gamma (a_{p+1}+l)} \over {\Gamma (b_1+l)\cdots \Gamma (b_p+l)\Gamma (1+l)}}}
=g_0(0)
\label{e16}
\end{equation}
\[-\sum\limits_{n=0}^N {g_n(s_p){1 \over {s_p+n}}{{(-s_p-n)_m} \over {\Gamma (m)}}}+T_1+T_2,\]
where
\begin{equation}
T_1=-\sum\limits_{l=m}^{\infty} {{{\Gamma (a_1+l)\Gamma (a_2+l)\cdots \Gamma (a_{p+1}+l)} \over {\Gamma (b_1+l)\cdots \Gamma (b_p+l)\Gamma (1+l)}}}
\label{e17}
\end{equation}
\[+\sum\limits_{n=0}^N {g_n(s_p)\sum\limits_{l=m}^\infty  {{{(-s_p-n)_l} \over {\Gamma (1+l)}}}},\]
\begin{equation}
T_2=-\mathop {\lim }\limits_{m\to \infty }\sum\limits_{n=N+1}^\infty  {{{g_n(s_p)} \over {(s_p+n)\Gamma (-s_p-n)}}{{\Gamma (-s_p-n+m)} \over {\Gamma (m)}}}.
\label{e18}
\end{equation}
We now have to show that the contributions $T_1$ and $T_2$ of the tails of the series are small. The terms in the series  in $T_2$ are $O(m^{-s_p-n})$ with $n\ge N+1$ as $m \to \infty$,  and so the series, in the worst case when $n=N+1$, converges provided that $\Re(-s_p-N)<0$. The sum of the series then is $O(m^{-s_p-N})$, and $T_2$ exists and vanishes. The term $T_1$ may be written
\begin{equation}
T_1=-\sum\limits_{l=m}^\infty  {{{\Gamma (-s_p+l)} \over {\Gamma (1+l)}}S_l},
\label{e19}
\end{equation}
where
\begin{equation}
S_l= {{{\Gamma (a_1+l)\Gamma (a_2+l)\cdots \Gamma (a_{p+1}+l)} \over {\Gamma (b_1+l)\cdots \Gamma (b_p+l)\Gamma (-s_p+l)}}}
\label{e20}
\end{equation}
\[-\sum\limits_{n=0}^N {{(-1)^n}g_n(s_p){1 \over {\Gamma (-s_p-n)(1+s_p-l)_n}}}.\]
In the second term, use has been made of the reflection formula of the gamma function.  It can be shown  \cite{b02} that $S_l$ is $O(l^{-N-1})$ as $l\to\infty$, so the terms of the sum are $O(l^{-s_p-N-2})$, and $T_1$ is $O(m^{-s_p-N-1})$ as $m\to\infty$. In this way we arrive at
\begin{theorem} 
If $s_p$ is not equal to an integer and  $N=0, 1, 2,...$ is chosen to be greater than $\Re(-s_p)$, then
\begin{equation}\sum\limits_{l=0}^{m-1} {{{\Gamma (a_1+l)\Gamma (a_2+l)\cdots \Gamma (a_{p+1}+l)} \over {\Gamma (b_1+l)\cdots \Gamma (b_p+l)\Gamma (1+l)}}}
=g_0(0)
\label{e21}
\end{equation}
\[-\sum\limits_{n=0}^N {(-1)^n{1 \over {s_p+n}}{{ (a_1+s_p)_n(a_2+s_p)_n} \over {(1)_n}}{{\Gamma (-s_p-n+m)} \over {\Gamma (m)}}}\]
\[\times \sum\limits_{k=0}^n {{{(-n)_k} \over {(a_1+s_p)_k(a_2+s_p)_k}}A_k^{(p)}}+O(m^{-s_p-N-1})\]
as $m \to\infty$.
\label{th3}
\end{theorem}
This theorem is applicable in the case of $\Re(s_p)>0$ too and then may give more details of the asymptotic behaviour than theorem \ref{th2}.

\section{The case when $s_p$ is equal to an integer}

If $s_p$ is equal to an integer, then the procedure of proof remains essentially the same, except that the starting point is a more complicated continuation formula, in place of (\ref{e12}), containing logarithmic terms which have to be expanded too. 

\subsection{The case when $s_p$ is equal to zero}
If $s_p$ is equal to zero, the hypergeometric series is called zero-balanced. It is this case which seems to be most interesting. While a few early results can be found in the monographs \cite{bai} \cite{sla} \cite{luk}, the topic had received renewed attention, with emphasis on the zero-balanced series,  in connection with Ramanujan's notebooks \cite{ram} \cite{be1} \cite{be2} \cite{est} \cite{eva} \cite{b87} \cite{sri} \cite{sai} \cite{sas} \cite{b98}.

The required continuation formula for the zero-balanced hypergeometric function is \cite{b92}
\begin{equation}
{{\Gamma (a_1)\Gamma (a_2)\cdots \Gamma (a_{p+1})} \over {\Gamma (b_1)\cdots \Gamma (b_p)}}{}_{p+1}F_p\left( \left.\begin{array}{c}a_1,a_2,\ldots ,a_{p+1}\\
  b_1,\ldots,b_p\end{array} \right| z \right)
\label{e22}
\end{equation} 
\[=\sum\limits_{n=0}^\infty  {d_n(1-z)^n}+\sum\limits_{n=0}^\infty  {e_n(1-z)^n\ln (1-z)},\]
where
\begin{equation}
e_n=-{{{(a_1)_n(a_2)_n} \over {(1)_n(1)_n}}\sum\limits_{k=0}^n {{{(-n)_k} \over {(a_1)_k(a_2)_k}}A_k^{(p)}}},
\label{e23}
\end{equation} 
\begin{equation}
d_0=2\psi (1)-\psi (a_1)-\psi (a_2)+\sum\limits_{k=1}^\infty  {{{\Gamma (k)} \over {(a_1)_k(a_2)_k}}A_k^{(p)}},
\label{e24}
\end{equation} 
and the other coefficients are not needed here. Again, we have to expand the singular terms on the right-hand side,
\begin{equation}
(1-z)^n\ln (1-z)=\sum\limits_{l=1}^\infty  {c_l^{(n)} z^l},
\label{e25}
\end{equation}
where
\begin{equation}
c_l^{(n)}=-{1 \over l}(-1)^n{{\Gamma(1+n)\Gamma (l-n)} \over {\Gamma (l)}}
\label{e26}
\end{equation}
for $l>n$, while the coefficients for $l\le n$ are 
\begin{equation}
c_1^{(1)}=-1,\quad c_1^{(2)}=-1,\quad c_2^{(2)}={\textstyle{3 \over 2}}.
\label{e27}
\end{equation}
With $z=1$, we need the partial sums of these series,
\begin{equation}
\sum\limits_{l=1}^{m-1} {c_l^{(0)}}=\psi (1)-\psi (m),
\label{e28}
\end{equation}
\begin{equation}
\sum\limits_{l=1}^{m-1} {c_l^{(1)}}=-{1 \over {m-1}},
\label{e29}
\end{equation}
\begin{equation}
\sum\limits_{l=1}^{m-1} {c_l^{(2)}}={1 \over {(m-1)(m-2)}}.
\label{e30}
\end{equation}
Keeping the contribution from the singular terms up to and including $n=2$ together with the constant term on the right-hand side of (\ref{e22}) we arrive at
\begin{theorem} 
If $s_p$ is equal to zero, then
\begin{equation}
\sum\limits_{l=0}^{m-1} {{{\Gamma (a_1+l)\Gamma (a_2+l)\cdots \Gamma (a_{p+1}+l)} \over {\Gamma (b_1+l)\cdots \Gamma (b_p+l)\Gamma (1+l)}}}=\sum\limits_{k=1}^\infty  {{{\Gamma(k) } \over {(a_1)_k(a_2)_k}}A_k^{(p)}}
\label{e31}
\end{equation}
\[+\psi (1)-\psi (a_1)-\psi (a_2)+\psi (m)+[a_1a_2-A_1^{(p)}](m-1)^{-1}\]
\[-{\textstyle{1 \over 4}}[a_1(a_1+1)a_2(a_2+1)-2(a_1+1)(a_2+1)A_1^{(p)}+2A_2^{(p)}][(m-1)(m-2)]^{-1}\]
\[+O(m^{-3})\]
as $m \to\infty$, where the infinite sum over $k$ converges if 
 $\Re(a_j)>0$ for $j=3, \ldots , p+1$. 
\label{th4}
\end{theorem}
This theorem gives more details of the asymptotic behaviour than our earlier formula \cite{b98}. If desired,  equation (\ref{e31}) may be rewritten using
\begin{equation}
\psi (m)=\ln (m)-{\textstyle{1 \over 2}}m^{-1}-{\textstyle{1 \over {12}}}m^{-2}+O(m^{-3})
\label{e32}
\end{equation}
in order to exhibit the logarithmic dependence on $m$.

\subsection{The case when $s_p$ is equal to a positive integer}

For $s_p$ equal to a positive integer $t$, the theorems \ref{th1} or \ref{th2} are applicable, but for  $s_p$ equal to $1$ or $2$ a few more terms of the asymptotic expansion might be desirable. In these cases the required continuation formula reads 
\cite{b92}
\begin{equation}
{{\Gamma (a_1)\Gamma (a_2)\cdots \Gamma (a_{p+1})} \over {\Gamma (b_1)\cdots \Gamma (b_p)}}{}_{p+1}F_p\left( \left.\begin{array}{c}a_1,a_2,\ldots ,a_{p+1}\\
  b_1,\ldots,b_p\end{array} \right| z \right)
\label{e33}
\end{equation} 
\[=\sum\limits_{n=0}^{t-1} {l_n(1-z)^n}+(1-z)^t\sum\limits_{n=0}^\infty  {[w_n+q_n\ln (1-z)](1-z)^n},\]
where
\begin{equation}
q_n=-(-1)^t{{(a_1+t)_n(a_2+t)_n} \over {\Gamma (1+t+n)}}\sum\limits_{k=0}^n {{{(-n)_k} \over {(a_1+t)_k(a_2+t)_k}}A_k^{(p)}},
\label{e34}
\end{equation} 
\begin{equation}
l_0={{\Gamma (a_1)\Gamma (a_2)\Gamma (t)} \over {\Gamma (a_1+t)\Gamma (a_2+t)}}\sum\limits_{k=0}^\infty  {{{(t)_k} \over {(a_1+t)_k(a_2+t)_k}}A_k^{(p)}},
\label{e35}
\end{equation} 
and the other coefficients are not needed here.
Using the expansions of the logarithmic terms and their partial sums from above and keeping the contributions up to and including $n=1$ in case of $t=1$ or $n=0$ in case of $t=2$, we may arrive at the following theorems. 
\begin{theorem} \label{th5}
If $s_p$ is equal to $1$, then
\begin{equation}
\sum\limits_{l=0}^{m-1} {{{\Gamma (a_1+l)\Gamma (a_2+l)\cdots \Gamma (a_{p+1}+l)} \over {\Gamma (b_1+l)\cdots \Gamma (b_p+l)\Gamma (1+l)}}}
\label{e36}
\end{equation}
\[={1 \over {a_1a_2}}\sum\limits_{k=0}^\infty  {{{(1)_k} \over {(a_1+1)_k(a_2+1)_k}}A_k^{(p)}}-(m-1)^{-1}\]
\[+{\textstyle{1 \over 2}}[(a_1+1)(a_2+1)-A_1^{(p)}][(m-1)(m-2)]^{-1}+O(m^{-3})\]
as $m \to\infty$,
where the infinite sum converges if 
$\Re(a_j)>0$ for $j=3, \ldots , p+1$. 
 \end{theorem}
\begin{theorem} \label{th6}
If $s_p$ is equal to $2$, then
\begin{equation}
\sum\limits_{l=0}^{m-1} {{{\Gamma (a_1+l)\Gamma (a_2+l)\cdots \Gamma (a_{p+1}+l)} \over {\Gamma (b_1+l)\cdots \Gamma (b_p+l)\Gamma (1+l)}}}
\label{e37}
\end{equation}
\[={1 \over {(a_1)_2(a_2)_2}}\sum\limits_{k=0}^\infty  {{{(2)_k} \over {(a_1+2)_k(a_2+2)_k}}A_k^{(p)}}-{\textstyle{1 \over 2}}[(m-1)(m-2)]^{-1}+O(m^{-3})\]
as $m \to\infty$, where the infinite sum converges if 
 $\Re(a_j)>0$ for $j=3, \ldots , p+1$. 
\end{theorem}
We may observe that theorem \ref{th5} and theorem \ref{th6} are not unexpected, but (\ref{e36}) and (\ref{e37}) are special cases of (\ref{e21}). Indeed, all the quantities in  (\ref{e21}) remain well-defined when $s_p$ approaches a positive integer, and so theorem \ref{th3} is valid even for positive integer values of $s_p$.

\subsection{The case when $s_p$ is equal to a negative integer}

When $s_p$ is equal to a negative integer $-t$, then the required continuation formula is
\cite{b92}
\begin{equation}
{{\Gamma (a_1)\Gamma (a_2)\cdots \Gamma (a_{p+1})} \over {\Gamma (b_1)\cdots \Gamma (b_p)}}{}_{p+1}F_p\left( \left.\begin{array}{c}a_1,a_2,\ldots ,a_{p+1}\\
  b_1,\ldots,b_p\end{array} \right| z \right)
\label{e38}
\end{equation} 
\[=(1-z)^{-t}\sum\limits_{n=0}^{t-1} {h_n(1-z)^n}+\sum\limits_{n=0}^\infty  {[u_n+v_n\ln (1-z)](1-z)^n},\]
where
\begin{equation}
h_n=(-1)^n{{(a_1-t)_n(a_2-t)_n\Gamma (t-n)} \over {\Gamma (1+n)}}\sum\limits_{k=0}^n {{{(-n)_k} \over {(a_1-t)_k(a_2-t)_k}}A_k^{(p)}},
\label{e39}
\end{equation}
\begin{equation}
v_n=-(-1)^t{{(a_1-t)_{t+n}(a_2-t)_{t+n}} \over {\Gamma (1+n)\Gamma (1+t+n)}}\sum\limits_{k=0}^{t+n} {{{(-t-n)_k} \over {(a_1-t)_k(a_2-t)_k}}A_k^{(p)}},
\label{e40}
\end{equation}
\begin{equation}
u_0=\sum\limits_{k=1}^\infty  {{{\Gamma (k)} \over {(a_1)_k(a_2)_k}}A_{t+k}^{(p)}}
\label{e41}
\end{equation}
\[+(-1)^t{{(a_1-t)_t(a_2-t)_t} \over {\Gamma (1+t)}}\sum\limits_{k=0}^{t} {{{(-t)_k} \over {(a_1-t)_k(a_2-t)_k}}A_k^{(p)}}\]
\[\times [\psi (1+t-k)+\psi (1)-\psi (a_1)-\psi (a_2)],\]
and the other $u_n$ are not needed here. Proceeding as above and keeping the contributions from the logarithmic terms up to and including $n=2$, we then may get 
\begin{theorem} \label{th7}
If $s_p$ is equal to a negative integer $-t$, then
\begin{equation}
\sum\limits_{l=0}^{m-1} {{{\Gamma (a_1+l)\Gamma (a_2+l)\cdots \Gamma (a_{p+1}+l)} \over {\Gamma (b_1+l)\cdots \Gamma (b_p+l)\Gamma (1+l)}}}
\label{e42}
\end{equation}
\[=\sum\limits_{n=0}^{t-1} {(-1)^n{1 \over {t-n}}{{(a_1-t)_n(a_2-t)_n\Gamma (t-n)} \over {\Gamma (1+n)}}{{\Gamma (t-n+m)} \over {\Gamma (m)}}}\]
\[\times\sum\limits_{k=0}^n {{{(-n)_k} \over {(a_1-t)_k(a_2-t)_k}}A_k^{(p)}}\]
\[+\sum\limits_{k=1}^\infty  {{{\Gamma (k)} \over {(a_1)_k(a_2)_k}}A_{t+k}^{(p)}}\]

\[+(-1)^t{{(a_1-t)_t(a_2-t)_t} \over {\Gamma (1+t)}}\sum\limits_{k=0}^t {{{(-t)_k} \over {(a_1-t)_k(a_2-t)_k}}A_k^{(p)}}\]
\[\times [\psi (1+t-k)-\psi (a_1)-\psi (a_2)+\psi (m)]\]
\[+(-1)^t{{(a_1-t)_{t+1}(a_2-t)_{t+1}} \over {\Gamma (2+t)}}\sum\limits_{k=0}^{t+1} {{{(-t-1)_k} \over {(a_1-t)_k(a_2-t)_k}}A_k^{(p)}}(m-1)^{-1}\]
\[-(-1)^t{{(a_1-t)_{t+2}(a_2-t)_{t+2}} \over {2\Gamma (3+t)}}\]
\[\times\sum\limits_{k=0}^{t+2} {{{(-t-2)_k} \over {(a_1-t)_k(a_2-t)_k}}A_k^{(p)}[}(m-1)(m-2)]^{-1}+O(m^{-3})\]
as $m \to\infty$, where the infinite sum converges if 
$\Re(a_j)>0$ for $j=3, \ldots , p+1$. 
\end{theorem}

\section {Special examples}

If the parameters $a_1, ..., a_{p+1}, b_1, ..., b_p$ satisfy certain additional relations or even have special values, it may well be that the appropriate formula above can be simplified, in particular its constant term, which involves an infinite sum.
 
As our first  example of this kind, let us consider a zero-balanced ${}_5F_4$ with the special values of the parameters
$a_1=a_2=a_3=a_4=1/2$, $a_5=5/4$, $b_1=b_2=b_3=1$, $b_4=1/4$. The constant term $CT$ of (\ref{e31}) then can be reduced (most conveniently by means of (\ref{e50})) to 
\begin{equation}
CT=-\psi (1)+4\ln (2)-{1 \over {18}}\sum\limits_{l=0}^\infty  {(1+l){{(2)_l(2)_l} \over {({\raise3pt\hbox{$\scriptstyle 5$} \!\mathord{\left/ {\vphantom {\scriptstyle {5 2}}}\right.\kern-\nulldelimiterspace} \!\lower3pt\hbox{$\scriptstyle 2$}})_l({\raise3pt\hbox{$\scriptstyle 5$} \!\mathord{\left/ {\vphantom {\scriptstyle {5 2}}}\right.\kern-\nulldelimiterspace} \!\lower3pt\hbox{$\scriptstyle 2$}})_l}}
{}_3F_2\left( \left.\begin{array}{c}
{\raise3pt\hbox{$\scriptstyle 3$} \!\mathord{\left/ {\vphantom {\scriptstyle {3 2}}}\right.\kern-\nulldelimiterspace} \!\lower3pt\hbox{$\scriptstyle 2$}},{\raise3pt\hbox{$\scriptstyle 3$} \!\mathord{\left/ {\vphantom {\scriptstyle {3 2}}}\right.\kern-\nulldelimiterspace} \!\lower3pt\hbox{$\scriptstyle 2$}},

-l\\ 2,3\end{array} \right| 1 \right) }.
\end{equation}
Here the infinite sum seems to be difficult to evaluate analytically, although it must simply be equal to $18\ln(2)$. For we know from Berndt \cite{be2} that
\begin{equation}
CT=-\psi (1)+3\ln (2).
\end{equation}
Using this information, we obtain from (\ref{e31})
\begin{corollary}
\label{th45}
\begin{equation}
{\textstyle{1 \over 4}}\pi ^2\sum\limits_{l=0}^{m-1} {{{({\raise3pt\hbox{$\scriptstyle 1$} \!\mathord{\left/ {\vphantom {\scriptstyle {1 2}}}\right.\kern-\nulldelimiterspace} \!\lower3pt\hbox{$\scriptstyle 2$}})_l({\raise3pt\hbox{$\scriptstyle 1$} \!\mathord{\left/ {\vphantom {\scriptstyle {1 2}}}\right.\kern-\nulldelimiterspace} \!\lower3pt\hbox{$\scriptstyle 2$}})_l({\raise3pt\hbox{$\scriptstyle 1$} \!\mathord{\left/ {\vphantom {\scriptstyle {1 2}}}\right.\kern-\nulldelimiterspace} \!\lower3pt\hbox{$\scriptstyle 2$}})_l({\raise3pt\hbox{$\scriptstyle 1$} \!\mathord{\left/ {\vphantom {\scriptstyle {1 2}}}\right.\kern-\nulldelimiterspace} \!\lower3pt\hbox{$\scriptstyle 2$}})_l({\raise3pt\hbox{$\scriptstyle 5$} \!\mathord{\left/ {\vphantom {\scriptstyle {5 4}}}\right.\kern-\nulldelimiterspace} \!\lower3pt\hbox{$\scriptstyle 4$}})_l} \over {(1)_l(1)_l(1)_l(1)_l({\raise3pt\hbox{$\scriptstyle 1$} \!\mathord{\left/ {\vphantom {\scriptstyle {1 4}}}\right.\kern-\nulldelimiterspace} \!\lower3pt\hbox{$\scriptstyle 4$}})_l}}}=\psi (m)-\psi (1)+3\ln (2)
\end{equation}
\[+{\textstyle{1 \over 4}}(m-1)^{-1}-{\textstyle{1 \over 8}}[(m-1)(m-2)]^{-2}+O(m^{-3}).\]
\end{corollary}

Our second example is a zero-balanced ${}_5F_4$ which depends only on three independent parameters $a$, $b$, $c$ according to
$a_1=a$, $a_2=b$, $a_3=c$, $a_4=a+b+c-1$, $a_5=(a+b+c+1)/2$, $b_1=b+c$, $b_2=c+a$, $b_3=a+b$, $b_4=(a+b+c-1)/2$. Here again, the infinite sum in the constant term of (\ref{e31}) is difficult to evaluate analytically, but fortunately we know the constant term from a different source \cite{eva}\cite{be2}\cite{ram}. Thus we obtain 
\begin{corollary}
\begin{equation}
\sum\limits_{l=0}^{m-1} {{{\Gamma (a+l)\Gamma (b+l)\Gamma (c+l)\Gamma (a+b+c-1+l)\Gamma ({\textstyle{1 \over 2}}[a+b+c+1]+l)} \over {\Gamma (b+c+l)\Gamma (c+a+l)\Gamma (a+b+l)\Gamma ({\textstyle{1 \over 2}}[a+b+c-1]+l)\Gamma (1+l)}}}
\end{equation}
\[=\psi (m)
+{\textstyle{1 \over 2}}[\psi (1)-\psi (a)-\psi (b)-\psi (c)]+{\textstyle{1 \over 2}}(a+b+c-1)(m-1)^{-1}
\]
\[
+{\textstyle{1 \over 4}}[2abc-(a+b+c)(a+b+c-1)][(m-1)(m-2)]^{-2}+O(m^{-3}).
\]
\end{corollary}
As a special case of this formula, with $a=b=c=1/2$, we can get again corollary \ref{th45} above.

Other partial sums with interrelated parameters have been obtained by A. K. Srivastava \cite{sri}.

\end{document}